\documentclass{amsart}

\date{20 May 2001}
\newtheorem{theorem}{Theorem}[section]
\newtheorem{prop}[theorem]{Proposition}
\newtheorem{lemma}[theorem]{Lemma}
\newtheorem{corollary}[theorem]{Corollary}

\def\supp{\text{supp }}
\def\C{\mathbb{C}}
\def\A{\mathbb{A}}
\def\D{\mathbb{D}}
\def\R{\mathbb{R}}

\def\S{\mathcal{S}}
\begin{document}

\title{The Banach envelope of Paley-Wiener type spaces}
\author{Mark Hoffmann}
\address{Department of Mathematics \\
University of Missouri-Columbia \\
Columbia, MO 65211\\ USA}
\email{mathgr26@math.missouri.edu}
\thanks{The author was partially supported by NSF grant DMS-9870027}

\begin{abstract}
We give an explicit computation of the Banach envelope for
 the Paley-Wiener type spaces $E^p,\, 0<p<1$.
 This answers a question by Joel Shapiro.
\end{abstract}

\maketitle

\section{Introduction}
The Paley-Wiener type space $E^p$
(a precise definition is given below)
 consists of certain band-limited functions \cite{pb}.
We will show that for $0<p<1$ this space can be identified as
a complemented subspace of the direct sum of two classical Hardy spaces.
Since the Banach envelopes of the Hardy spaces are known, we are able to
establish a necessary and sufficient condition for entire functions to belong
to the envelope $E^p_c$.

For an open subset $\Omega\subseteq \C$ let $\A(\Omega)$ be the space
of holomorphic functions on $\Omega$.
An entire function $f$ is of exponential type $\tau > 0$ if for all
 $\epsilon > 0$ there is a
$C_\epsilon >0$ such that for all $z\in \C$ we have
$|f(z)|\leq C_\epsilon e^{(\tau+\epsilon) |z|}$.
For $0<p<\infty$ let $E^{p,\tau}$ be the space of entire functions of
exponential type $\tau$ such that their restrictions to the real axis
 are in $L^p(\R) $ :
 $$E^{p,\tau}=\left\{ f \in \A(\C) : f \text{ has exponential type $\tau$ },
 f|_{\R}\in L^p(\R)\right\}.$$
 We will only consider $\tau=\pi$ and write $E^p=E^{p,\pi}$ from now on.
 This causes no loss of generality, because
  we can simply rescale a function $f\in E^{p,\tau}$ to obtain
   $\supp (f|_\R)^\wedge \subseteq [-\pi,\pi]$.
   The quantity
  $$\|f\|_{E^p}=\left(\int_{-\infty}^{\infty}|f(x)|^p dx\right)^{1/p}
  =\|f|_{\R}\|_{L^p}$$
  defines a norm on $E^p$ for $1\leq p<\infty $ and a quasi-norm for $0<p<1$.

 These spaces $E^p$ are complete, and hence
 can be identified with closed subspaces of $L^p(\R)$ (e.g. see
 \cite{ce} or \cite{tr}). A classical
 theorem of Paley and Wiener gives a characterization of $E^2$ as the image
  of the inverse Fourier transform of $L^2[-\pi,\pi]$. Hence
   functions in
  $E^2$ have compactly supported Fourier transform, i.e. they are
  band-limited.
 These functions are important to signal processing due to
   their sampling properties (Shannon sampling theorem).
   Since
  $\|f\|_q\leq C_{pq}\|f\|_p$ for $f\in E^p, 0<q\leq p$ \cite{tr}, we have
  $E^p\subset E^q$ for $0<p\leq q$. In particular,
   $E^p\subset E^1\subset E^2$ for $0<p<1$.
  A characterization of $E^p$ for $0<p<1$ as a discrete Hardy
  space is shown in \cite{ce}.

Let $X$ be a quasi-normed space with separating dual .
Then the Banach envelope $X_{c}$
of $X$
 is the completion of $(X,\|\cdot\|_{C})$ where $C$ is
 the convex hull of the closed unit ball $B_X$, and $\|\cdot \|_{C}$
  is the
 Minkowski functional of $C$. $X_{c}$
 is a Banach space. The Banach envelope
  is characterized up to isomorphism by
   $(X_c)^*=X$ and $\overline{X}=X_c$. Every operator $T:X\rightarrow Y$
   extends uniquely to $\tilde{T}:X_c\rightarrow Y_c$
(\cite{kp}, \cite{mm}, \cite{js} and \cite{wo}).

The standard example for finding a Banach envelope is $\ell^p$ for $0<p<1$.
In this case we have $\ell^p\subset \ell^1$, $\ell^p$ is dense in $\ell^1$ and
$(\ell^p)^*=\ell^\infty=(\ell^1)^*$. Therefore, $l^p_c=l^1$. The spaces $E^p$ are nested as well,
$E^p\subset E^1$ for $0<p<1$. Furthermore, $E^p$ is dense in $E^1$ since it contains
all Schwartz functions with Fourier transform supported in $[-\pi,\pi]$. This makes
$E^1$ a candidate for the envelope of $E^p$, but
it turns out that $E^p_c$ is a certain weighted Bergman space of entire
functions different from $E^1$. The proof relies heavily on the thoery of Hardy
spaces $H^p$. We use a deep result by Duren, Romberg and Shields \cite{pr}
 that the
Banach envelope of $H^p$ over the unit disk is a certain weighted $L^1$-Bergman
space (see below).

  The problem of identifying the
  Banach envelope of $E^p$ as a space of entire functions was originally
  posed by Joel Shapiro.
  I would like to thank Professor Nigel Kalton for communicating this problem and his
  helpful suggestions.
  I would also like to thank the referee for some very constructive comments.

\section{Preliminaries}
We recall the classical Hardy spaces of the disc $\D=\{z\in \C :
|z|\leq 1\}$
and the upper half planes $\C_\pm=\{z=x+iy\in \C : y\in \R_\pm\}$
$$H^p(\D)=\left\{f\in \A(\D) :
\|f\|^p_{H^p(\D)}=\sup_{0\leq r<1}\int_{0}^{2\pi} |f(re^{i\theta})|^p
\frac{d\theta}{2\pi}< \infty\right\},$$
$$H^p(\C_+)=\left\{f\in \A(\C_+) :
\|f\|^p_{H^p(\C_+)}=\sup_{y>0}\int_{-\infty}^\infty |f(x+iy)|^p dx
< \infty\right\}.$$
Analogously we define for the lower half plane
$$H^p(\C_-)=\left\{f\in \A(\C_-) :\|f\|^p_{H^p(\C_-)}=
\sup_{y<0}\int_{-\infty}^\infty |f(x+iy)|^p dx < \infty\right\}.$$
An isometric isomorphism between $H^p(\C_+)$ and $H^p(\C_-)$ is given
by $f(z)\mapsto \overline{f(\bar z)}$.
Let $\S$ be the space of Schwartz functions and
 $\S'$ the space of tempered distributions.
Every $f\in H^p(\C_+)$ is uniquely determined by its
boundary value distribution
$f^b=\lim_{y\rightarrow 0} f(x+iy)\in \S'.$
Denote the space of these boundary distributions by $H^p_+(\R)$, and let
$\|f^b\|_{H^p_+(\R)}=\|f\|_{H^p(\C_+)}$.
In the same way define $H^p_-(\R)$. The real
and imaginary part of $f$ have a boundary value distribution in
the real Hardy space as defined in \cite{cw}, \cite{st}.
Hence the Fourier transform of $f\in H^p_+(\R), 0<p<1$, is a
 continuous functions and satisfies
$\supp \widehat{f} \subseteq [0,\infty)$. More precisely, we have
$|\widehat{f}(\xi)|\leq C|\xi|^{1/p-1}\|f\|_{H^p_+(\R)}$. Transfer to the
lower half plane shows that
$f\in H^p_-(\R)$ has $\supp\hat{f}\subseteq (-\infty,0]$.
All these results can be found e.g. in
\cite{cw}, \cite{gc}, \cite{ko}, \cite{st}.\\

The Bergman spaces over the disc
and the upper half plane are defined for $0<p<\infty, \, \alpha > -1$ as
$$A^{p,\alpha}(\D)=\left\{f\in\A(\D) : \|f\|_{p,\alpha}^p=
\int_{\D}|f(x+iy)|^p(1-|z|^2)^\alpha
dx\,dy
<\infty\right\},$$
$$A^{p,\alpha}(\C_+)=\left\{f\in\A(\C_+) :  \|f\|_{p,\alpha}^p=
\int_{\C_+}|f(x+iy)|^p y^\alpha dx\,dy<\infty\right\}.$$
Analogously we define $A^{p,\alpha}(\C_-).$
The Banach envelope of $H^p(\D)$ was identified by Duren, Romberg
and Shields \cite{pr}.
 \begin{prop} $H^p_c(\D)=A^{1_,1/p-2}(\D).$
 \end{prop}
For a different approach in the setting of Besov and Triebel-Lizorkin spaces
 spaces see \cite{mm} and \cite{mi}. In particular, it is shown that
for a bounded Lipschitz domain $\Omega\subset \R^2$
and $1/2 < p < 1$ we have
\begin{prop}
$H^p_c(\Omega)=A^{1,1/p-2}(\Omega).$
\end{prop}
We will need the analogous statement for the upper and lower half plane.
The following proposition is certainly well-known, for completeness
we give a proof using a conformal map from $\C_+$ onto $\D$.
\begin{prop} \label{res2}
$H^p_{c}(\C_\pm)=A^{1,1/p-2}(\C_\pm)$
\end{prop}
\begin{proof}
It is enough to consider the upper half plane.
We use the conformal map $w=\phi(z)=\frac{i-z}{i+z}$ from $\C_+$ onto
$\D.$ Then with $F(z)=f(\phi(z)), z=x+iy$ we have
$f(w)\in H^p(\D)$ if and only if $F(z)/(z+i)^{2/p}\in H^p(\C_+)$
(see \cite{pd}, \cite{ko}).
For $\alpha > -1$ we get
$$\int_\D|f(w)|^p (1-|w|^2)^\alpha dx\,dy
=\int_{\C_+}|F(z)|^p (1-|\phi(z)|^2)^\alpha |\phi'(z)|^2 dx\,dy.$$
A short computation shows
$ 1-|\phi(z)|^2=\frac{4y}{|z+i|^2}$ and $|\phi'(z)|^2=\frac{4}{|z+i|^4}$.
Therefore,
$$\int_\D|f(w)|^p (1-|w|^2)^\alpha dx\,dy
=4^{\alpha+1}\int_{\C_+}\left|\frac{F(z)}{(z+i)^{(2\alpha+4)/p}}\right|^p
y^\alpha dx\,dy.$$
Hence $F(z)\in H^p_{c}(\C_+)$ if and only if
$f(w)(\phi^{-1}(w)+i)^{2/p}\in A^{1,\alpha}(\D)$ where $\alpha=1/p-2$.
Then $2\alpha+4=2/p$, and mapping back to the upper half plane shows
$F(z)\in H^p_{c}(\C_+)$ if and only if $F(z)\in A^{1,\alpha}(\C_+)$.
\end{proof}

\section{The Banach envelope of $E^p$}
The next proposition is due to Plancherel and Polya \cite{pp}.
\begin{prop} \label{planp}
Let $0<p<\infty$ and $f\in E^p$. Then for every $y\in\R$ we have
$$\int_{-\infty}^\infty |f(x+iy)|^p dx \leq e^{p\pi |y|}
\int_{-\infty}^\infty |f(x)|^p dx.$$
\end{prop}
For an entire function $f$ let $f_{\pm\pi}=e^{\pm i\pi z}f|_{\C_\pm}(z)$
 and $j(f)=(f_{\pi},f_{-\pi})$.
Then from \ref{planp} it follows that
$f\mapsto f_{\pm\pi}$ is an isometric isomorphism of $E^p$ into $H^p(\C_{\pm}).$
Hence $j$ embeds $E^p$ into $H^p(\C_+)\oplus H^p(\C_-).$\\

Crucial to compute the envelope of $E^p$ is the following.
\begin{lemma}
$j(E^p)$ is complemented in $H^p(\C_+)\oplus H^p(\C_-).$
\end{lemma}
\noindent
\begin{proof}
We construct a bounded projection onto $j(E^p).$
Choose $\phi,\psi\in \S$ such that
$\supp \hat{\phi}\subset [-2\pi,\pi], \supp \hat{\psi}\subset [-\pi,2\pi]$
 and
$\hat{\phi}(x)+\hat{\psi}(x)=1$ for all $x\in [-\pi,\pi]$. This can be done
by a suitable partition of unity on the Fourier transform side. Then
let $T: \S'\times \S'\mapsto \S'$ be defined by
 $$T(u,v)=u*\phi + v*\psi.$$
We have $\supp T(u,v)^\wedge\subseteq [-2\pi,2\pi]$,
and if $u$ has $\supp \hat{u}\subseteq [-\pi,\pi]$ then $T(u,u)=u.$
Write $u_{\pm\pi}=e^{\pm i\pi x}u$ for $u\in \S'$.
Suppose $(u,v)\in H^p_+({\R})\oplus H^p_-({\R})$. Then
$T(u_{-\pi},v_\pi)^\wedge$ is continuous and supported in $[-\pi,\pi]$.
This shows that $T(u_{-\pi},v_\pi)$ has an extension to a function in
 $E^2$. Non-tangential (distributional) boundary values of functions
 in $H^p(\C_+)$ are in $L^p(\R)$,
 and we have $\|u*\Phi\|_{H^p_+(\R)}\leq C\|u\|_{H^p_+(\R)}$
  for $u\in H^p_+(\R), \Phi\in \S'$
 \cite{st}. Hence $T(u_{-\pi},v_\pi)|_\R \in L^p(\R)$, and
 $T(u_{-\pi},v_\pi)$ extends to a function in $E^p$.
 This extension has the explicit form
 $Q : H^p(\C_+)\oplus H^p({\C_-}) \mapsto E^p$,
$$Q(f,g)(z)=\left< T(f_{-\pi}^b,g_{\pi}^b), e^{itz} \right>
=\int_{-\pi}^\pi T(f_{-\pi}^b,g_{\pi}^b)e^{itz}dt$$
By choice of $\phi,\psi$ we have $Q(j(f))=f$,
and hence $P=jQ$ is the desired projection.
\end{proof}

Now we arrive at our characterization of $E^p_{c,q}$.
\begin{theorem} \label{res} An entire function $f$ belongs to
 $E^p_{c}$  if and only if
$$ \|f\|=\int_{\C} e^{-\pi |y|}|y|^{1/p-2} |f(x+iy)| dx\,dy
<\infty .$$
Moreover, $\|\cdot\|$ is equivalent to the quasi-norm of $E^p_{c}$.
\end{theorem}
\noindent
\begin{proof}   Let $\alpha=1/p-2$.
Define
$$Z=\{f\in \A(\C) : f_{\pi}\in A^{1,\alpha}({\C_+}),
 f_{-\pi}\in A^{1,\alpha}({\C_-})\}\subset A^{1,\alpha}(\C_+)\oplus
 A^{1,\alpha}(\C_-).$$
We will show $E^p_{c}=Z$ with equivalence of norms.
We have $H^p(\C_\pm)\subset A^{1,\alpha}(\C_\pm)$, and hence $E^p\subset Z$.
It is crucial to observe that the operator $Q$ from the previous proof
extends
to $\tilde{Q}: A^{1,\alpha}(\C_+)\oplus A^{1,\alpha}(\C_-)\mapsto E^p_{c}$
while preserving the defining equation, i.e.
$\tilde{Q}(f,g)(z)=\langle T(f^b_{-\pi},g^b_\pi),e^{itz} \rangle$.
This follows from the characterization of boundary
value distributions for functions in $A^{1,\alpha}({\C_+})$ \cite{rt}.
These distributions are uniquely determined by their values on $\S$, and
passing from $f\in A^{q,\alpha}({\C_+})$ to $f^b\in A^{1,\alpha}_+({\R})$
 is a continuous operation.
Therefore, we have $E^p_{c}\subset E^2$ and
$\tilde{Q}(j(f))=f$ for $f\in Z$.
We conclude that $j(E^p_{c})$ is a complemented subspace of
$H^p_{c}(\C_+)\oplus H^p_{c}(\C_-)=A^{1,\alpha}(\C_+)\oplus
A^{1,\alpha}(\C_-)$.
 Hence $E^p_{c}\subseteq Z$, and on $j(E^p_{c})$ the envelope-norm is
  equivalent to the norm of $A^{1,\alpha}(\C_+)\oplus A^{1,\alpha}(\C_-)$.
  It remains to establish density of $j(E^p)$ in $j(Z)$.
Pick $f\in Z$. Then there are $(h_n,g_n)\in
 H^p(\C_+)\oplus H^p(\C_-)$ such that $(h_n,g_n)\rightarrow
 (f_{\pi},f_{-\pi})=j(f)\in A^{1,\alpha}(\C_+)\oplus A^{1,\alpha}(\C_-)$.
 Then $Q(h_n,g_n)\in E^p$ and $Q(h_n,g_n)\rightarrow
 \tilde{Q}(j(f))=f$.
\end{proof}
Let us see that $E^p_c\neq E^1$.
For $f\in E^1, y>0$ let $f_y(x)=f(x+iy)$. Then by \ref{planp} we have
$\|f_y\|_{L^1}\leq e^{\pi y}\|f_0\|_{L^1}=\|f\|_{E^1}$. For a given $y>0$ we
choose a function that satisfies the converse inequality up to a constant as
follows: Take $\phi_\epsilon\in \S$ such that
$\supp \phi_\epsilon\subseteq [-\pi,-\pi+\epsilon]$, and let $f^\epsilon(z)
=\left<\phi_\epsilon,e^{itz}\right>$.
Then for fixed $y_0>0$ we choose $\epsilon>0$ small enough to give
$\|f^\epsilon_y\|_{L^1}\geq C e^{\pi y}\|f^\epsilon\|_{E^1}$ for all
 $0\leq y\leq y_0$.
We obtain
$\|f^\epsilon\|_{E^p_c}=
\int_{-\infty}^\infty e^{-\pi |y|}|y|^{1/p-2}\|f^\epsilon_y\|_{L^1}dy
\geq C\int_{-\infty}^{y_0} y^{1/p-2}\|f^\epsilon\|_{E^1}dy
= (1/p-1)^{-1}y_0^{1/p-1}\|f^\epsilon\|$.
Since $1/p-1>0$ we conclude that the $E^1$-norm and the $E^p_c$-norm are not
equivalent.

A result from \cite{kt} (see also \cite{wo})
is $A^{1,\alpha}(\D)\approx \ell^1$,
 and therefore we have using a result of
 Pe\l czy\'nski \cite{Pe} that
 every complemented subspace of $\ell^1$ is isomorphic to $\ell^1$:
\begin{corollary} $E^p_{c}\approx \ell^1.$
\end{corollary}
\section{The $q$-envelope of $E^p$}
The result from the previous section can be generalized to the $q$-envelopes.
 Let $X$ be a quasi-normed space with separating dual and $0<q\leq 1$.
Then the Banach q-envelope $X_{c,q}$
of $X$
 is the completion of $(X,\|\cdot\|_{C_q})$ where $C_q$ is
 the $q-$convex hull of the closed unit ball $B_X$, and $\|\cdot \|_{C_q}$
  is the
 Minkowski functional of $C_q$. $X_{c,q}$ is a complete
 quasi-normed space.

The results in \cite{cr} and \cite{wo} give the $q$-envelopes of the Hardy spaces.
\begin{prop}
For $0<p<q\leq 1$ we have
$H^p_{c,q}(\D)=A^{q,q/p-2}(\D).$
\end{prop}

The proofs of \ref{res} and \ref{res2} work analogously for the $q$-envelope, $0<q<1$,
 if we use $\alpha=q/p-2$ and $A^{q,\alpha}(\C_\pm)$. Hence
\begin{theorem} An entire function $f$ belongs to
 $E^p_{c,q}$  if and only if
$$ \|f\|=\left(\int_{\C} e^{-q\pi |y|}|y|^{q/p-2} |f(x+iy)|^q dx\,dy
\right)^{1/q}<\infty .$$
Moreover, $\|\cdot\|$ is equivalent to the quasi-norm of $E^p_{c,q}$.
\end{theorem}
In \cite{kt} (see also \cite{wo})
it is shown that $A^{q,\alpha}(\D)\approx \ell^q$,
 and therefore we have using a result of Stiles \cite{stil} or
 \cite{kp} for $q<1$ that every complemented subspace of $\ell^q$ is isomorphic to $\ell^q$:
\begin{corollary} $E^p_{c,q}\approx \ell^q.$
\end{corollary}


\begin{thebibliography}{99}
\bibitem {pb} R. P. Boas, \it Entire functions\rm, Academic Press,
 New York, 1954.
\bibitem{cw} R. R. Coifman and R. Rochberg, \it
Representation theorems for holomorphic and harmonic functions in $L_p$\rm,
Ast\'erisque \bf 77 \rm (1980) 11-66.
\bibitem{cr} R. R. Coifman and G. Weiss, \it Extensions of Hardy spaces
 and their use in analysis\rm, Bull. Amer. Math. Soc. \bf 83 \rm (1977) 569-645.
\bibitem{pd} P. L. Duren, {\it Theory of $H_p$ spaces}, Academic Press,
New York/London, 1970.
\bibitem{pr} P. L. Duren, B. W. Romberg and A. L. Shields ,
\it Linear functionals on $H_p$ spaces when $0<p<1$\rm, J. Reine Angew.
Math. {\bf 238} (1969) 32-60.
\bibitem{ce} C. Eoff, \it The discrete nature of the Paley-Wiener spaces\rm,
  Proc. Amer. Math. Soc. \bf 123 \rm (1995) 505-512.
\bibitem{gc} J. Garcia-Cuerva and J. L. Rubio de Francia, {\it Weighted Norm
 Inequalities and Related Topics}, North-Holland, 1985.
\bibitem{kp} N. J. Kalton, N. T. Peck and J. W. Roberts,
{\it An F-space Sampler}, London Mathematical Society, LNS No. \bf 89 \rm, 1984.
\bibitem{kt} N. J. Kalton and D. A. Trautman, {\it Remarks on
subspaces of $H_p$ when $0<p<1$}, Michigan Math. J. {\bf 29} (1982)
163-171.
\bibitem{ko} P. Koosis, {\it Introduction to $H_p$ spaces} ,
Cambridge Univ. Press, 2nd ed. 1998.
\bibitem{mm} O. Mendez and M. Mitrea, \it The Banach envelopes
of Besov and Triebel-Lizorkin spaces and applications to partial differential
equations\rm, J. Fourier Anal. Appl. {\bf 6} (2000) 503-533.
\bibitem{mi} M. Mitrea, \it Banach envelopes of holomorphic Hardy spaces\rm,
preprint, 2000.
\bibitem{Pe} A. Pe\l czy\'nski, \it Projections in certain Banach spaces, \rm Studia Math. {\bf 19} (1960) 209--228.


\bibitem{pp} M. Plancherel and G. P\'olya, \it Fonctions enti\`eres et
int\`egrales de Fourier multiples\rm, Comment. Math. Helv. \bf 10 \rm (1937)
110-163.
\bibitem{rt} F. Ricci and M. Taibleson, \it Boundary Values of Harmonic
Functions in Mixed Norm Spaces and Their Atomic Structure\rm,
Ann.  Scuola Norm. Sup. Pisa \bf 10 \rm (1983) 1-54.
\bibitem{tr} H. Triebel, {\it Theory of function spaces II}, Birkh\"auser,
Berlin, 1992.
\bibitem{js} J. H. Shapiro, \it Mackey topologies, reproducing kernels, and
diagonal maps on the Hardy and Bergman spaces\rm, Duke Math. J. \bf 43 \rm
 (1976) 187-202.
\bibitem{st} E. M. Stein, {\it Harmonic analysis: Real-Variable Methods,
Orthogonality, and Oscillatory integrals}, Princeton Univ. Press,
 Princeton, NJ, 1993.
\bibitem{stil} W. J. Stiles, \it Some properties of $l\sb{p}$, $0<p<1$, \rm Studia Math. {\bf 42} (1972) 109-119. 
\bibitem{wo} P. Wojtaszczyk, \it $H_p$-spaces, $p\leq 1$, and spline
systems\rm , Studia Math. {\bf 77} (1984) 289-320.
\end{thebibliography}
\end{document}